\documentclass{amsart}
\usepackage{amsfonts,amsmath,amssymb,eucal,tipa}
\usepackage[all]{xy}
\usepackage{hyperref}
\usepackage[utf8]{inputenc}

\begin{document}

\parindent=0pt
\parskip=6pt

\newcommand{\A}{{\mathbb A}}
\newcommand{\G}{{\mathbb G}}
\newcommand{\Pb}{{\mathbb P}^1({\mathbb R})}
\newcommand{\Gl}{{\rm Gl}}
\newcommand{\F}{{\mathbb F}}
\newcommand{\C}{{\mathbb C}}
\newcommand{\E}{{\mathbb E}}
\newcommand{\I}{{\mathbb I}}
\newcommand{\pP}{{\mathbb P}}
\newcommand{\Ll}{{\mathbb L}}
\newcommand{\Q}{{\mathbb Q}}
\newcommand{\R}{{\mathbb R}}
\newcommand{\T}{{\mathbb T}}
\newcommand{\Z}{{\mathbb Z}}
\newcommand{\Sl}{{\rm Sl}}
\newcommand{\PGl}{{\rm PGl}}
\newcommand{\fc}{{\mathfrak c}}

\newcommand{\SU}{{\rm SU}}
\newcommand{\pc}{{+_\circ}}
\newcommand{\sslash}{{/\!\!/}}  

\newcommand{\ie}{\textit{ie}\,}
\newcommand{\eg}{\textit{eg}\,}
\newcommand{\se}{{\sf e}}
\newcommand{\cf}{{\textit{cf}\,}}

\newcommand{\cay}{{\mathfrak c}}
\newcommand{\bR}{{\overline{\mathbb R}}}
\newcommand{\sgn}{{\rm sgn}}
\newcommand{\slg}{{\rm slog}}
\newcommand{\maps}{{\rm maps}}
\newcommand{\Aff}{{{\mathbb A}{\rm ff}}}
\newcommand{\wG}{{{\widetilde{\mathbb G}}}}
\newcommand{\cR}{{\mathcal R}}
\newcommand{\tw}{{\widetilde{\tau}}}
\newcommand{\bc}{{\bf c}}


\title{On the universal cover of the real projective line}

\author[Jack Morava]{Jack Morava}

\address{Department of Mathematics, The Johns Hopkins University,
Baltimore, Maryland}

\begin{abstract}{Retrieved from the lost libraries of Atlantis}\end{abstract}

\maketitle 

{\bf \S I} The composition law
\[
x,y \mapsto x +_\pP y = \frac{x + y}{1 - xy} \;,
\]
defines an abelian group structure on the Pythagorean line
\[
 \pP := \pP^1(\R) = (\R^2 - \{0\})/\R^\times \cong \R \cup \infty
 \]
with $0$ as unit and elements\begin{footnote}{in terms of the cross-ratio $[x_0:x_1:x_2:x_3] := \frac{x_{01}}{x_{02}} \frac{x_{23}}{x_{13}} : \overline{\mathfrak M}_{0,4}(\R) \cong \Pb$}\end{footnote}
\[
x = [x:1] = [0:x:1:\infty] 
\]
eg\begin{footnote}{$x \mapsto x +_\pP  \infty$ is continuous but not smooth at $x=0$}\end{footnote} $x +_\pP  \infty = - x^{-1}$; its completion at 0 is a one-dimensional formal group law over $\Q$ with $x \mapsto \arctan x$ as its logarithm and $x +_\pP (-x) = 0$. Its underlying topological space is the circle $[-\infty,+\infty]/\{\pm \infty\}$ with the endpoints consolidated as a single element of order two.\bigskip

{\sc Proposition} {\bf 1} {\it There is a commutative diagram
\[
\xymatrix{
0 \ar[r] & \Z \ar[d]^{-\cay} \ar[r] & \Z \times_\pm \R^\times \ar[d] \ar@{.>}[r] & \Pb \ar[d]^= \ar[r] & 0 \\
0 \ar[r] & \pP^\vee \ar[r] & \wG_m(\R) \ar[r]^\tw & \pP \ar[r] & 0 }
\]
of locally compact abelian groups and continuous homomorphisms, as defined below.}\bigskip


We need three preliminary definitions: 

$\vdash i$) $\bR = \R \cup {\pm \infty}$ is the end-point compactification of the reals, while $\pP = \R \cup \infty$ is its Alexandroff (one-point) compactification. We regard $\tan : \R \to \bR$ as $\pi$-periodic, and 
we take
\[ 
{\rm Arctan} : \bR \to [-\pi/2, +\pi/2] : x \mapsto \frac{i}{2} \log \frac{1 - ix}{1 + ix}  \in \Q \{ \{x\} \} 
\]
to be the analytic function defined by its power series on $\R \subset \bR$. Then 
\[
[\tan \circ {\rm Arctan} : \bR \to [-\pi/2,+\pi/2] \subset \R \to \bR] := {\bf 1}_\bR
\]
is the identity, while 
\[
[{\rm Arctan} \circ \tan : \R \to \bR \to [-\pi/2,+\pi/2] \subset \R] := {\bf 1}_{[-\pi/2,+\pi/2]} 
\]
is the periodic (sawtooth) extension of the identity map of the interval $[-\pi/2,+\pi/2]$. 

$\vdash ii$) For any real $q >0$ let 
\[
[Q] := [q^{-1/2},q^{+1/2}] \subset \R
\]
denote the closed interval of width $q^{-1/2}(q-1)$; note that the involution $x \mapsto r(x) = x^{-1}$ takes this interval to itself. The map
\[
\cR^\times := \{\pm 1\} \times \cR_{>0} = \{\pm 1\} \times \{q^\Z\} \cdot [Q] \to  \wG_m(\R) \to \{\pm 1\} \times q^\Z [Q]/\langle s \rangle \cong \G_m(\R) 
\]
(defined by real multiplication) is equivariant with respect to the action of the infinite dihedral group $I_\infty := \langle r, s \: | \: rsr = s^{-1}  \rangle$ (defined by $s(x) = qx$ , $r(x) = x^{-1}$) on both its domain and range. 
 
If $q = \exp(\pi)$, then
\[
\tau(x) := \tan \log |x| :  [q^{-1/2},q^{+1/2}]  = [Q] \to \R \cup \infty = \pP
\]
is monotone and surjective. It satisfies $\tau(s(x)) = \tau(qx) = \tau(x)$, defining a quotient covering
\[
\tw:  \{\pm 1\} \times \cR_{>0} / \langle s \rangle  \to  \wG_m(\R) \to \{\pm 1\} \times (q^Z \cdot [Q])/\langle s \rangle \to \pP \; .
\]

$\vdash iii$) The Cayley transform \cite{6}
\[
[\cay](x) =  \left[\begin{array}{cc}
                             1 & -i \\
                             1 & +i 
                                    \end{array}\right](x) = \frac{x-i}{x+i}  : \bR \to \T
\]
(\ie stereographic projection) satisfies 
\[
[\cay](x) =  -\exp(2i \cdot \frac{i}{2} \log \frac{1 - ix}{1 + ix} )  =   - \exp (2i {\rm Arctan} \; x) .
\]
The matrix $\cay \in M_2(\Z[i])$ has determinant $2i$, trace $2i^{1/2}$, eigenvalues $2$ and $2i$; it represents an 
element of order four in ${\rm PGl}_2(\C)$; alternatively, 
\[
 \bc :=  \left[\begin{array}{cc}
                                              1 & -1 \\
                                               1 & 1
                                    \end{array}\right] \in \Sl_2(\Q(\surd 2))
\]
(\ie $[\cay](ix) = [{\bc}](x)$) exchanges $\{\pm 1\}$ and $\{0,\infty\}$: $[\bc](-x) = [\bc](x)^{-1}$, $[\bc](e^x) = - [\bc](e^{-x})  \dots$). \bigskip

{\sc Proof:} of the proposition:  The assertions follow by direct calculation:
\[
\tau(xy) = \tau(x) +_\pP \tau(y) \;,
\]
while
\[
 [\cay](x +_\pP y) = - [\cay](x) \cdot [\cay](y).
\]

 [{\sc corollary}  {\it Pontrjagin's dual character group
\[
\pP^\vee = {\rm Hom}_c (\pP, \T) \cong \Z \cdot \cay
\]
($\T \subset \C^\times$) of the real projective line is generated by $-\cay$.}

For $\exp(2it \log |x|) : \R^\times \to \T$ factors through $\pP$ only if $t \in \pi \Z$.]
$\Box$ 

{\sc Proposition} {\bf 2} {\it The diagram
\[
\xymatrix{ 
0 \ar[r] & \Z \ar[d] \ar[r] & \wG_m(\R) \ar[d]^\log \ar@{.>}[dr]^\chi \ar[r]^\tw & \pP \ar[d]^{-\cay} \ar[r] &  0 \\
0 \ar[r] & \Z_2 \ar[r] & \R^\times/q^\Z \ar[r]^-2 & \T \cong \R/2\pi i \Z \ar[r] & 0 }
\]
(of locally compact groups and continuous homomorphisms) commutes}. 

Here
\[
\chi(x) = ([\cay] \circ \tw)(x) = \exp(2i \log |x|) = |x|^{2i} : 
\R^\times_{>0}/q^\Z  \to \R/{2\pi i}\Z  \cong \T  \subset \C^\times 
\]
is an unramified Tate quasicharacter of $\R^\times$.

{\sc exercise} \cite{2}, \cite{5}(\S 2,3) The relativistic Mercator logarithm
\[
\lambda(x) := {\rm Arctanh} \circ \sin(x) = - \log i [\cay] (e^{ix}) = - \log i [\bc]( -i \exp (ix)) \dots 
\]\bigskip

{\sc definition} Let us call the group of fractional linear transformations 
\[
\R \times \R^\times \ni (s,t) : x  \mapsto  
\left[\begin{array}{cc}
                  s & t \\
                  0 & 1
                 \end{array}\right](x) = sx + t \;,
\]      
the affine group $\Aff(\R)$ of $\R$. The analogous subgroup  $\Aff(\Z) \cong I_\infty$ is dihedral.\bigskip

{\sc Proposition} {\bf 3} {\it The diagram
\[
\xymatrix{
0 \ar[r] & \R \ar[d]^\Z \ar@{.>}[dr] \ar[r] & \Aff(\R) \cong \R \rtimes \G_m(\R)  \ar[d]^\Z \ar[r] & \G_m(\R) \ar[d]^{\mu_2} \ar[r] & 0 \\ 
0 \ar[r] & \T \ar[r] & {\mathbb E} := \Aff(\R)/\Aff(\Z) \ar[r] & \pP \ar[r] & 0 },
\]
defines a morphism of fibrations, if not of groups. The sequence 
\[
\xymatrix{
0 \ar[r] & \pi_1 \T \ar[d] \ar[r] & \pi_1\E \ar@{.>}[d] \ar[r] & \pi_1 \pP \ar[d] \ar[r] & 0 \\
0 \ar[r] & \Z \ar[r] & \Z \rtimes \Z  \ar[r] & \Z \ar[r] & 0 }
\]
of fundamental groups identifies $\pi_1\E$ as the semidirect product defined by the sign $\Z \to \mu_2$.}

 $\E$ represents a generator of $H^1(\pP,\Z) = \pi_0 \maps(\pP,B\Z) \cong \Z$ as the mapping torus of the antipodal map of the circle, which unwraps as an infinite sequence of Dehn twists.
 
{\sc Corollary} {\it The total space of the circle bundle $\E$ over the Pythagorean line $\pP$ is a Klein surface \dots }\bigskip

{\sc ? Corollary IIUC} {\it The universal cover 
\[
\C \to \C/\Z[i] \to \E =\Aff(\R)/\Aff(\Z) : r e^{2\pi it} \mapsto [z \mapsto rz + t] 
\]
factors through the elliptic curve with $j$-invariant 1728} \bigskip

{\sc Remarks:}  This proposition presents $\E$ as a circle bundle (analogous to periodic time) over a one-dimensional Pythagorean/Archimedean space/topos \cite{10}. Abel-Jacobi theory \cite{8} identifies its orientation cover with itself, now carrying a level two structure
\[
H_1(\C/\Z[i],\Z) \equiv \Z^2 
\]
(\ie modulo two) which assigns $\E$ a spin structure \cite{1}, \cf further \cite{7}. The splitting 
\[
t \ni \pi_1\T \cong \Z  \to (t,s) \ni \pi_1\C/\Z[i] \cong \Z^2  \to \pi_1\Pb \ni s
\]
presents the action of the mapping-class group $\pi_0{\rm Diff}_\pm (\C/\Z[i]) \cong \Gl_2(\Z)$ on that Jacobian. See further \cite{3,4} \bigskip

 {\bf \S II  The Two Toruses}\bigskip
      
{\bf 2.1} If we write
\[
x \mapsto [J(x)] := 
\left[\begin{array}{cc}
                   1 & x \\
                 - x & 1
           \end{array}\right] \in {\rm PGl}_2(\R) 
\]
then 
\[
[J(x +_\pP y)] =
\left[\begin{array}{cc}
                   1 & \frac{x + y}{1 - xy}  \\
                 - \frac{x + y}{1 - xy}  & 1
            \end{array}\right]  =  [J(x)] \circ [J(y)]
\]
even though 
\[
\frac{\det J(x)J(y)}{\det J(x+_\pP y)}  = (1 - xy)^2 ,
\]
so 
\[ 
[J(\tau(xy))] = [J(\tau(x) +_\pP \tau(y))] = [J(\tau(x)] \circ [J(\tau(y)] 
\]
defines a homomorphism 
\[
\pP_1(\R) \ni x \mapsto [J(\tau(x))] \in {\rm PGl}_2(\R) .
\]

{\bf 2.2} The map
\[
\C \ni X \mapsto H(X) = \left[\begin{array}{cc}
                                                    1 +iX & 0 \\
                                                     0 & 1 -iX 
                                                     \end{array}\right] \in M_2(\C)[x]
\]
satisfies $H(X) \circ H(Y)  = (1 -XY) H(X +_\pP Y)$,
defining something like a hyperbolic torus $\C^\times \to \PGl_2(\C)$. Unfortunately it blows up at $X = \pm i  \in \C$.\bigskip

The identification (beware determinants) of 
\[
\left[\begin{array}{cc}
                  \alpha - \delta & -\beta + \gamma \\
                  \beta + \gamma  & \alpha + \delta
        \end{array}\right] \in \Sl_2(\R) 
\] 
with 
\[        
\left[\begin{array}{cc}
                  \alpha + i \beta & \gamma + i \delta  \\
                  \gamma - i \delta  & \alpha - i \beta
        \end{array}\right] \in \SU(1,1)
\]
sends 
\[(1 + x^2)^{-1/2} \cdot   
\left[\begin{array}{cc}
                   1 & x \\
                  - x  & 1
        \end{array}\right]  \mapsto (1 + x^2)^{-1/2} \cdot       
\left[\begin{array}{cc}
                  1 + ix & 0  \\
                  0 & 1 - ix
        \end{array}\right] ,
 \]
ie
\[
\mapsto 
\left[\begin{array}{cc}
                  e^{i\theta} & 0  \\
                  0 & e^{-i\theta}
        \end{array}\right] ,
\]
when $x = \tan \theta$. 

{\bf 2.3 Claim} {\it In the diagram
\[
\xymatrix{
\R^\times \ar[d]^-\log \ar[r]^-\tau & \Pb \times 0 \ar[d] \ar[r]^-J & \SU(1,1)/\{\pm 1\} \ar@{.>}[d] \\
\C^\times \cong \R \times \T \ar[r]^-{\tan \times 1_\T} & \Pb \times \T \ar[r]^-{H \times *_\T} & \PGl_2(\C) }
\]
the representation $\C(-1) \oplus \C(1)$ of $\C^\times$ defined by the homomorphism on the lower side of the diagram is the complexification of the representation of the torus $\R^\times$ defined by the composition along the upper side.} \bigskip

{\bf Remarks} 

{\bf 1} Following Lebesgue, an increasing bijection $u : [0,1] \to [0,1]$ is differentiable almost everywhere, defining a compact topological abelian group law
\[
+_u : x,y \to x +_u y = u^{-1}(u(x) + u(y))
\]
on $\R/\Z$, possibly of interest in probability theory, 

{\bf 2} The proposition implies that the fundamental complex representation of $\SU(1,-1)(\R) \cong \Sl_2(\R)$ pulls back to the Pythagorean torus $\Pb$. If we let $V= \{1, u, v, w \}$ be Klein's four-group\begin{footnote}{$u^2 = v^2 = \dots =1, uv = w, vw = u \dots$}\end{footnote}, with $u$ interpreted as the nontrivial central element $-1 \in \Sl_2(\R)$, then 
\[
\SU(1,-1)(\C) \ltimes_{\pm 1} V \to {\rm O}(1,3) 
\]
is a spin cover of the full Lorentz group, making $V$ analogous to the group $\{1,P,T,PT\}$ generated by space and time reversal : to be distinguished from Wu's group (generated by parity and charge reversal) of symmetries in elementary particle physics \cite{9,10}.

\end{document}